# SPIKE SOLUTIONS FOR A CLASS OF SINGULARLY PERTURBED QUASILINEAR ELLIPTIC EQUATIONS

MARCO SQUASSINA

ABSTRACT. By means of a penalization scheme due to del Pino and Felmer, we prove the existence of single–peaked solutions for a class of singularly perturbed quasilinear elliptic equations associated with functionals which lack of smoothness. We don't require neither uniqueness assumptions on the limiting autonomous equation nor monotonicity conditions on the nonlinearity. Compared with the semilinear case some difficulties arise and the study of concentration of the solutions needs a somewhat involved analysis in which the Pucci–Serrin variational identity plays an important role.

## 1. Introduction and the main result

Let $\Omega$ be a possibly unbounded smooth domain of $\mathbb{R}^N$ with $N \geqslant 3$.

Since the pioneering work of Floer and Weinstein [FW] in the one space dimension, much interest has been directed in the last decade to singularly perturbed elliptic problems of the form

$$
(1) \quad \begin{cases} -\varepsilon^2 \Delta u + V(x)u = f(u) & \text{in } \Omega \\ u > 0 & \text{in } \Omega \\ u = 0 & \text{on } \partial\Omega \end{cases}
$$

for a superlinear and subcritical nonlinearity $f$ with $f(s)/s$ nondecreasing.

Typically, there exists a family of solutions $(u_\varepsilon)_{\varepsilon>0}$ which exhibits a spike shape around the local minima (possibly degenerate) of the function $V(x)$ and decade elsewhere as $\varepsilon$ goes to zero (see e.g. [ABC, DF, DF2, DF3, JT2, O, O2, R, SP, SP2, W] and references therein). A natural question is now whether these concentration phenomena are a special feature of the semilinear case or we can expect a similar behaviour to hold for more general elliptic equations which possess a variational structure.

Date: October 25, 2018.
1991 *Mathematics Subject Classification.* 35J40; 58E05.
This work was partially supported by Ministero dell'Università e della Ricerca Scientifica e Tecnologica (40% – 2001) and by Gruppo Nazionale per l'Analisi Funzionale e le sue Applicazioni.

In this paper we will give a positive answer to this question for the following class of singularly perturbed quasilinear elliptic problems

$$
(2) \quad \begin{cases} -\varepsilon^2 \sum_{i,j=1}^{N} D_j(a_{ij}(x,u)D_i u) + \frac{\varepsilon^2}{2} \sum_{i,j=1}^{N} D_s a_{ij}(x,u)D_i u D_j u + V(x)u = f(u) & \text{in } \Omega \\ u > 0 & \text{in } \Omega \\ u = 0 & \text{on } \partial\Omega \end{cases}
$$

under suitable assumptions on the functions $a_{ij}$, $V$ and $f$. Notice that if $a_{ij}(x,s) = \delta_{ij}$ then equation (2) reduces to (1), in which case the problem originates from different physical and biological models and, in particular, in the study of the so called *standing waves* for the nonlinear Schrödinger equation.

Existence and multiplicity results for equations like (2) have been object of a very careful analysis since 1994 (see e.g. [AB, AB2, C, CD, S2] for the case where $\Omega$ is bounded and [CG, S] for $\Omega$ unbounded). On the other hand, to the author's knowledge, no result on the *asymptotic behaviour* of the solutions (as $\varepsilon$ vanishes) of (2) can be found in literature. In particular no achievement is known so far concerning the concentration phenomena for the solutions $u_\varepsilon$ of (2) around the local minima, not necessarily nondegenerate, of $V$.

We stress that various difficulties arise in comparison with the study of the semilinear equation (1) (see Section 5 for a list of properties which are not known to hold in our framework).

A crucial step in proving our main result is to show that the *Mountain–Pass* energy level of the functional $J$ associated with the autonomous limiting equation

$$
(3) \quad -\sum_{i,j=1}^{N} D_j(a_{ij}(\widehat{x},u)D_i u) + \frac{1}{2}\sum_{i,j=1}^{N} D_s a_{ij}(\widehat{x},u) D_i u D_j u + V(\widehat{x}) u = f(u) \quad \text{in } \mathbb{R}^N
$$

with $\widehat{x} \in \mathbb{R}^N$, is the *least* among other nontrivial critical values (Lemma 3.5). Notice that, no uniqueness result is available, to our knowledge, for this general equation (on the contrary in the semilinear case some uniqueness theorems for ground state solutions have been obtained by performing an ODE analysis in radial coordinates, see e.g. [CL]). The least energy problem for (3) is also related to the following fact:

$$
(4) \quad u \in H^1(\mathbb{R}^N),\ u \geqslant 0 \text{ and } u \text{ solution of } (3) \text{ implies that } J(u) = \max_{t \geqslant 0} J(tu)
$$

Unfortunately, as remarked in [CG, section 3], if one assumes that condition (10) holds, then property (4) *cannot* hold true even if the map $\{s \mapsto f(s)/s\}$ is nondecreasing.

In order to show the minimality property for the Mountain–Pass level and to study the uniform limit of $u_\varepsilon$ on $\partial\Lambda$, inspired by the recent work of Jeanjean and Tanaka [JT], we make a repeated use of the Pucci–Serrin identity [PS], which has turned out to be a very powerful tool (Lemmas 3.5 and 3.6).

Notice that the functional associated with (2) (see (15)) is not even locally Lipschitz and tools of nonsmooth critical point theory will be emploied (see [CDM, DM] and references therein). Also the proof of a suitable Palais–Smale type condition for a modification of the functional $I_\varepsilon$ becomes more involved.

We assume that $f \in C^1(\mathbb{R}^+)$ and there exist $1 < p < \frac{N+2}{N-2}$ and $2 < \vartheta \leqslant p+1$ with

(5) $$\lim_{s \to +\infty} \frac{f(s)}{s^p} = 0, \qquad \lim_{s \to 0^+} \frac{f(s)}{s} = 0,$$

(6) $$0 < \vartheta F(s) \leqslant f(s)s \quad \text{for every } s \in \mathbb{R}^+,$$

where $F(s) = \int_0^s f(t)\,dt$ for every $s \in \mathbb{R}^+$.

Furthermore, let $V : \mathbb{R}^N \to \mathbb{R}$ be a locally Hölder continuous function bounded below away from zero, that is, there exists $\alpha > 0$ with

(7) $$V(x) \geqslant \alpha \quad \text{for every } x \in \mathbb{R}^N.$$

The functions $a_{ij}(x,s) : \Omega \times \mathbb{R}^+ \to \mathbb{R}$ are continuous in $x$ and of class $C^1$ with respect to $s$, $a_{ij}(x,s) = a_{ji}(x,s)$ for every $i,j = 1, \ldots, N$ and there exists a positive constant $C$ with

$$|a_{ij}(x,s)| \leqslant C, \quad |D_s a_{ij}(x,s)| \leqslant C$$

for every $x \in \Omega$ and $s \in \mathbb{R}^+$. Finally, let $R, \nu > 0$ and $0 < \gamma < \vartheta - 2$ be such that

(8) $$\sum_{i,j=1}^N a_{ij}(x,s)\xi_i\xi_j \geqslant \nu|\xi|^2,$$

(9) $$\sum_{i,j=1}^N s D_s a_{ij}(x,s)\xi_i\xi_j \leqslant \gamma \sum_{i,j=1}^N a_{ij}(x,s)\xi_i\xi_j,$$

(10) $$s \geqslant R \implies \sum_{i,j=1}^N D_s a_{ij}(x,s)\xi_i\xi_j \geqslant 0$$

for every $x \in \Omega$, $s \in \mathbb{R}^+$ and $\xi \in \mathbb{R}^N$.

Hypothesis (5), (6) and (7) on $f$ and $V$ are standard. Observe that neither *monotonicity* assumptions on the function $f(s)/s$ nor *uniqueness* conditions on the limiting equation (3) are considered. Finally, (9) and (10) have already been used, for instance in [AB, AB2, C, CD, CG], in order to tackle these general equations.

Let $H_V(\Omega)$ be the weighted Hilbert space defined by

$$H_V(\Omega) = \left\{ u \in H_0^1(\Omega) : \int_\Omega V(x) u^2 < +\infty \right\},$$

endowed with the scalar product $(u,v)_V = \int_\Omega Du Dv + V(x)uv$ and denote by $\|\cdot\|_{H_V(\Omega)}$ the corresponding norm.

Let $\Lambda$ be a compact subset of $\Omega$ such that there exists $x_0 \in \Lambda$ with

(11) $$V(x_0) = \min_\Lambda V < \min_{\partial \Lambda} V,$$

(12) $$\sum_{i,j=1}^N a_{ij}(x_0, s)\xi_i\xi_j = \min_{x \in \Lambda} \sum_{i,j=1}^N a_{ij}(x,s)\xi_i\xi_j$$

for every $s \in \mathbb{R}^+$ and $\xi \in \mathbb{R}^N$. Moreover, let us set

$$\sigma := \sup \big\{ s > 0 : \ f(t) \leqslant tV(x_0) \text{ for every } t \in [0,s] \big\} \tag{13}$$

and

$$\mathscr{M} := \big\{ x \in \Lambda : \ V(x) = V(x_0) \big\}. \tag{14}$$

The following is the main result of the paper.

**Theorem 1.1.** *Assume that conditions* (5), (6), (7), (8), (9), (10), (11), (12) *hold.*
*Then there exists $\varepsilon_0 > 0$ such that, for every $\varepsilon \in (0, \varepsilon_0)$, there exist $u_\varepsilon \in H_V(\Omega) \cap C(\overline{\Omega})$ and $x_\varepsilon \in \Lambda$ satisfying the following properties:*

(a) $u_\varepsilon$ *is a weak solution of the problem*

$$(P_\varepsilon) \quad \begin{cases} -\varepsilon^2 \sum_{i,j=1}^N D_j(a_{ij}(x,u)D_i u) + \frac{\varepsilon^2}{2} \sum_{i,j=1}^N D_s a_{ij}(x,u) D_i u D_j u + V(x)u = f(u) & \text{in } \Omega \\ u > 0 & \text{in } \Omega \\ u = 0 & \text{on } \partial\Omega\,; \end{cases}$$

(b) *there exists $\sigma' > 0$ such that*

$$u_\varepsilon(x_\varepsilon) = \sup_\Omega u_\varepsilon, \quad \sigma < u_\varepsilon(x_\varepsilon) < \sigma', \quad \lim_{\varepsilon \to 0} d(x_\varepsilon, \mathscr{M}) = 0$$

*where $\sigma$ is as in* (13) *and $\mathscr{M}$ is as in* (14) *;*

(c) *for every $\varrho > 0$ we have*

$$\lim_{\varepsilon \to 0} \|u_\varepsilon\|_{L^\infty(\Omega \setminus B_\varrho(x_\varepsilon))} = 0\,;$$

(d) *we have*

$$\lim_{\varepsilon \to 0} \|u_\varepsilon\|_{H_V(\Omega)} = 0$$

*and, as a consequence,* $\lim_{\varepsilon \to 0} \|u_\varepsilon\|_{L^q(\Omega)} = 0$ *for every $2 \leqslant q < +\infty$.*

The proof of the theorem is variational and in the spirit of a well–known paper by del Pino and Felmer [DF], where it was succesfully developed into a *local* setting the global approach initiated by Rabinowitz [R].

We will consider the functional $I_\varepsilon : H_V(\Omega) \to \mathbb{R}$ associated with the problem $(P_\varepsilon)$,

$$I_\varepsilon(u) := \frac{\varepsilon^2}{2} \sum_{i,j=1}^N \int_\Omega a_{ij}(x,u) D_i u D_j u + \frac{1}{2} \int_\Omega V(x) u^2 - \int_\Omega F(u) \tag{15}$$

and construct a new functional $J_\varepsilon$ which satisfies the Palais–Smale condition (in a suitable sense) at every level ($I_\varepsilon$ does not, in general) and to which the (nonsmooth) Mountain–Pass Theorem can be directly applied to get a critical point $u_\varepsilon$ with precise energy estimates.

Then we will prove that $u_\varepsilon$ goes to zero uniformly on $\partial \Lambda$ as $\varepsilon$ goes to zero (this is the hardest step, here we repeatedly use the Pucci–Serrin identity in a suitable form) and show that $u_\varepsilon$ is actually a solution of the original problem with all of the stated properties.

**Remark 1.2.** We do not know whether the solutions of problem $(P_\varepsilon)$ obey to the following *exponential decay*

$$(16) \qquad u_\varepsilon(x) \leqslant \alpha \exp\left\{-\frac{\beta}{\varepsilon}|x - x_\varepsilon|\right\} \quad \text{for every } x \in \Omega, \text{ for some } \alpha, \beta \in \mathbb{R}^+,$$

which is a typical feature in the semilinear case. This fact would follow if we had a suitable Gidas–Ni–Nirenberg [GNN] type result for the equation (3) to be combined with some results by Rabier and Stuart [RS] on the exponential decay of second order elliptic equations.

**Remark 1.3.** As pointed out in [DF3], the concentration around the minima of the potential is, in some sense, a *model situation* for other phenomena such as concentration around the maxima of $d(x, \partial\Omega)$. Furthermore it seems to be the *technically simplest* case, thus suitable for a first investigation in the quasilinear case.

The organization of the paper is as follows:

– In Section 2 we construct the modified functional $J_\varepsilon$ and we prove that it satisfies a variant of the classical Palais–Smale condition (see Definition 6.6).

– In Section 3 we study the concentration of the solutions $u_\varepsilon$ around the local minimum of $V(x)$ as $\varepsilon$ tends to zero.

– In Section 4 we finally end up the proof of Theorem 1.1.

– In Section 5 we list a few open problems related to our main result.

– In Section 6 we quote from [CD] various tools of nonsmooth critical point theory.

## 2. The del Pino–Felmer penalization scheme

We now define a suitable modification of the functional $I_\varepsilon$ in order to regain the (concrete) Palais–Smale condition at any level and apply Proposition 6.9 of the Appendix for every $\varepsilon > 0$. Let us consider the positive constant

$$\ell := \sup\left\{s > 0 : \frac{f(t)}{t} \leqslant \frac{\alpha}{k} \text{ for every } 0 \leqslant t \leqslant s\right\}$$

for some $k > \vartheta/(\vartheta - 2)$. We define the function $\widetilde{f} : \mathbb{R}^+ \to \mathbb{R}$ by setting

$$\widetilde{f}(s) := \begin{cases} \frac{\alpha}{k}s & \text{if } s > \ell \\ f(s) & \text{if } 0 \leqslant s \leqslant \ell \end{cases}$$

and the maps $g, G : \Omega \times \mathbb{R}^+ \to \mathbb{R}$,

$$g(x, s) := \chi_\Lambda(x)f(s) + (1 - \chi_\Lambda(x))\widetilde{f}(s), \quad G(x, s) = \int_0^s g(x, \tau)\, d\tau$$

for every $x \in \Omega$. Then the function $g(x,s)$ is measurable in $x$, of class $C^1$ in $s$ and it satisfies the following assumptions:

$$\text{(17)} \qquad \lim_{s \to +\infty} \frac{g(x,s)}{s^p} = 0, \qquad \lim_{s \to 0^+} \frac{g(x,s)}{s} = 0 \quad \text{uniformly in } x \in \Omega,$$

$$\text{(18)} \qquad 0 < \vartheta G(x,s) \leqslant g(x,s)s \quad \text{for every } x \in \Lambda \text{ and } s \in \mathbb{R}^+,$$

$$\text{(19)} \qquad 0 \leqslant 2G(x,s) \leqslant g(x,s)s \leqslant \frac{1}{k}V(x)s^2 \quad \text{for every } x \in \Omega \setminus \Lambda \text{ and } s \in \mathbb{R}^+.$$

Without loss of generality, we may assume that

$$g(x,s) = 0 \quad \text{for every } x \in \Omega \text{ and } s < 0,$$
$$a_{ij}(x,s) = a_{ij}(x,0) \quad \text{for every } x \in \Omega,\ s < 0 \text{ and } i,j = 1,\ldots,N.$$

Let $J_\varepsilon : H_V(\Omega) \to \mathbb{R}$ be the functional

$$J_\varepsilon(u) := \frac{\varepsilon^2}{2} \sum_{i,j=1}^N \int_\Omega a_{ij}(x,u) D_i u D_j u + \frac{1}{2} \int_\Omega V(x) u^2 - \int_\Omega G(x,u).$$

The next result provides the link between the critical points of the modified functional $J_\varepsilon$ and the solutions of the original problem.

**Proposition 2.1.** *Assume that $u_\varepsilon \in H_V(\Omega)$ is a critical point of $J_\varepsilon$ and that there exists a positive number $\varepsilon_0$ such that*

$$u_\varepsilon(x) \leqslant \ell \quad \text{for every } \varepsilon \in (0,\varepsilon_0) \text{ and } x \in \Omega \setminus \Lambda.$$

*Then $u_\varepsilon$ is a solution of $(P_\varepsilon)$.*

*Proof.* By assertion $(a)$ of Corollary 6.8, it results that $u_\varepsilon$ is a solution of (62). Since $u_\varepsilon \leqslant \ell$ on $\Omega \setminus \Lambda$, we have

$$G(x, u_\varepsilon(x)) = F(u_\varepsilon(x)) \quad \text{for every } x \in \Omega.$$

Moreover, by arguing as in the proof of [S, Lemma 1], one gets $u_\varepsilon > 0$ in $\Omega$.

Then $u_\varepsilon$ is a solution of $(P_\varepsilon)$. $\square$

The next Lemma – which is nontrivial – provides a local compactness property for bounded concrete Palais–Smale sequences of $J_\varepsilon$ (see Definition 6.6). For the proof, we refer the reader to [S, Theorem 2 and Lemma 3].

**Lemma 2.2.** *Assume that conditions* (5), (6), (7), (8), (9), (10) *hold. Let $\varepsilon > 0$. Assume that $(u_h) \subset H^1(\mathbb{R}^N)$ is a bounded sequence and*

$$\langle w_h, \varphi \rangle = \varepsilon^2 \sum_{i,j=1}^N \int_{\mathbb{R}^N} a_{ij}(x,u_h) D_i u_h D_j \varphi + \frac{\varepsilon^2}{2} \sum_{i,j=1}^N \int_{\mathbb{R}^N} D_s a_{ij}(x,u_h) D_i u_h D_j u_h \varphi$$

*for every $\varphi \in C_c^\infty(\mathbb{R}^N)$, where $(w_h)$ is strongly convergent in $H^{-1}(\widetilde{\Omega})$ for a given bounded domain $\widetilde{\Omega}$ of $\mathbb{R}^N$.*

*Then $(u_h)$ admits a strongly convergent subsequence in $H^1(\widetilde{\Omega})$. In particular, if $(u_h)$ is a bounded concrete Palais–Smale condition for $J_\varepsilon$ at level $c$ and $u$ is its weak limit, then, up to a subsequence, $Du_h \to Du$ in $L^2(\widetilde{\Omega}, \mathbb{R}^N)$ for every bounded subset $\widetilde{\Omega}$ of $\Omega$.*

Since $\Omega$ may be unbounded, in general, the original functional $I_\varepsilon$ does not satisfy the concrete Palais–Smale condition. In the following Lemma we prove that, instead, the functional $J_\varepsilon$ satisfies it for every $\varepsilon > 0$ at every level $c \in \mathbb{R}$.

**Lemma 2.3.** *Assume that conditions (5), (6), (7), (8), (9), (10) hold. Let $\varepsilon > 0$.*
*Then $J_\varepsilon$ satisfies the concrete Palais–Smale condition at every level $c \in \mathbb{R}$.*

*Proof.* Let $(u_h) \subset H_V(\Omega)$ be a concrete Palais–Smale sequence for $J_\varepsilon$ at level $c$.
We divide the proof into two steps:

<u>STEP I</u>. Let us prove that $(u_h)$ is bounded in $H_V(\Omega)$. Since $J_\varepsilon(u_h) \to c$, from inequalities (18) and (19), we get

$$(20) \quad \frac{\vartheta \varepsilon^2}{2} \sum_{i,j=1}^N \int_\Omega a_{ij}(x, u_h) D_i u_h D_j u_h + \frac{\vartheta}{2} \int_\Omega V(x) u_h^2$$
$$\leqslant \int_\Lambda g(x, u_h) u_h + \frac{\vartheta}{2k} \int_{\Omega \setminus \Lambda} V(x) u_h^2 + \vartheta c + o(1)$$

as $h \to +\infty$. Moreover, in view of Proposition 6.4, we have $J'_\varepsilon(u_h)(u_h) = o(\|u_h\|_{H_V(\Omega)})$ as $h \to +\infty$. Then, again by virtue of (19), we deduce that,

$$\varepsilon^2 \sum_{i,j=1}^N \int_\Omega a_{ij}(x, u_h) D_i u_h D_j u_h + \frac{\varepsilon^2}{2} \sum_{i,j=1}^N \int_\Omega D_s a_{ij}(x, u_h) u_h D_i u_h D_j u_h$$
$$+ \int_\Omega V(x) u_h^2 \geqslant \int_\Lambda g(x, u_h) u_h + o(\|u_h\|_{H_V(\Omega)}),$$

as $h \to +\infty$, which, by (9), yields

$$(21) \quad \left(\frac{\gamma}{2} + 1\right) \varepsilon^2 \sum_{i,j=1}^N \int_\Omega a_{ij}(x, u_h) D_i u_h D_j u_h + \int_\Omega V(x) u_h^2$$
$$\geqslant \int_\Lambda g(x, u_h) u_h + o(\|u_h\|_{H_V(\Omega)})$$

as $h \to +\infty$. Then, in view of (8), by combining inequalities (20) and (21) one gets

$$(22) \quad \min\left\{\left(\frac{\vartheta}{2} - \frac{\gamma}{2} - 1\right) \nu \varepsilon^2, \frac{\vartheta}{2} - \frac{\vartheta}{2k} - 1\right\} \int_\Omega \left(|Du_h|^2 + V(x) u_h^2\right)$$
$$\leqslant \vartheta c + o(\|u_h\|_{H_V(\Omega)}) + o(1)$$

as $h \to +\infty$, which implies the boundedness of $(u_h)$ in $H_V(\Omega)$.

<u>STEP II</u>. By virtue of Step I, there exists $u \in H_V(\Omega)$ such that, up to a subsequence, $(u_h)$ weakly converges to $u$ in $H_V(\Omega)$.

Let us now prove that actually $(u_h)$ converges strongly to $u$ in $H_V(\Omega)$. By taking into account Lemma 2.2 (applied with $\widetilde{\Omega} = B_\varrho(0)$ for every $\varrho > 0$), it suffices to prove that for every $\delta > 0$ there exists $\varrho > 0$ such that

$$\text{(23)} \qquad \limsup_h \int_{\Omega \setminus B_\varrho(0)} \left( |Du_h|^2 + V(x)u_h^2 \right) < \delta.$$

We may assume that $\Lambda \subset B_{\varrho/2}(0)$. Consider a cut–off function $\psi_\varrho \in C^\infty(\Omega)$ with $\psi_\varrho = 0$ on $B_{\varrho/2}(0)$, $\psi_\varrho = 1$ on $\Omega \setminus B_\varrho(0)$, $|D\psi_\varrho| \leqslant c/\varrho$ on $\Omega$ for some positive constant $c$. Let $M$ be a positive number such that

$$\text{(24)} \qquad \left| \frac{1}{2} \sum_{i,j=1}^N D_s a_{ij}(x,s) \xi_i \xi_j \right| \leqslant M \sum_{i,j=1}^N a_{ij}(x,s) \xi_i \xi_j$$

for every $x \in \Omega$, $s \in \mathbb{R}^+$, $\xi \in \mathbb{R}^N$ and let $\zeta : \mathbb{R} \to \mathbb{R}$ be the map defined by

$$\text{(25)} \qquad \zeta(s) := \begin{cases} 0 & \text{if } s < 0 \\ Ms & \text{if } 0 \leqslant s < R \\ MR & \text{if } s \geqslant R, \end{cases}$$

being $R > 0$ the constant defined in (10). Notice that

$$\text{(26)} \qquad \sum_{i,j=1}^N \left[ \frac{1}{2} D_s a_{ij}(x,s) + \zeta'(s) a_{ij}(x,s) \right] \xi_i \xi_j \geqslant 0, \quad \text{for every } x \in \Omega, s \in \mathbb{R}, \xi \in \mathbb{R}^N.$$

By Proposition 6.4, we can compute $J_\varepsilon'(u_h)(\psi_\varrho u_h \exp\{\zeta(u_h)\})$. Since $(u_h)$ is bounded in $H_V(\Omega)$ and (26) holds, we get

$$o(1) = J_\varepsilon'(u_h)(\psi_\varrho u_h \exp\{\zeta(u_h)\})$$

$$= \varepsilon^2 \sum_{i,j=1}^N \int_\Omega a_{ij}(x, u_h) D_i u_h D_j u_h \psi_\varrho \exp\{\zeta(u_h)\}$$

$$+ \varepsilon^2 \sum_{i,j=1}^N \int_\Omega \left[ \frac{1}{2} D_s a_{ij}(x, u_h) + \zeta'(u_h) a_{ij}(x, u_h) \right] D_i u_h D_j u_h u_h \psi_\varrho \exp\{\zeta(u_h)\}$$

$$+ \varepsilon^2 \sum_{i,j=1}^N \int_\Omega a_{ij}(x, u_h) u_h D_i u_h D_j \psi_\varrho \exp\{\zeta(u_h)\} + \int_\Omega V(x) u_h^2 \psi_\varrho \exp\{\zeta(u_h)\}$$

$$- \int_\Omega g(x, u_h) u_h \psi_\varrho \exp\{\zeta(u_h)\} \geqslant \int_\Omega \left( \varepsilon^2 \nu |Du_h|^2 + V(x) u_h^2 \right) \psi_\varrho \exp\{\zeta(u_h)\}$$

$$+ \varepsilon^2 \sum_{i,j=1}^N \int_\Omega a_{ij}(x, u_h) u_h D_i u_h D_j \psi_\varrho \exp\{\zeta(u_h)\} - \int_\Omega g(x, u_h) u_h \psi_\varrho \exp\{\zeta(u_h)\}.$$

Therefore, in view of (19), it results

$$o(1) \geqslant \int_\Omega \left(\varepsilon^2 \nu |Du_h|^2 + V(x)u_h^2\right) \psi_\varrho \exp\{\zeta(u_h)\}$$
$$+ \varepsilon^2 \sum_{i,j=1}^N \int_\Omega a_{ij}(x,u_h) u_h D_i u_h D_j \psi_\varrho \exp\{\zeta(u_h)\} - \frac{1}{k} \int_\Omega V(x) u_h^2 \psi_\varrho \exp\{\zeta(u_h)\}$$

as $\varrho \to +\infty$. Taking into account that

$$\left| \sum_{i,j=1}^N \int_\Omega a_{ij}(x,u_h) u_h D_i u_h D_j \psi_\varrho \exp\{\zeta(u_h)\} \right| \leqslant \frac{\exp\{MR\}\widetilde{C}}{\varrho} \|Du_h\|_2 \|u_h\|_2,$$

there exists $C' > 0$ (which depends only on $\varepsilon, \nu$ and $k$) such that, as $\varrho \to +\infty$,

$$\limsup_h \int_{\Omega \setminus B_\varrho(0)} \left(|Du_h|^2 + V(x)u_h^2\right) \leqslant \frac{C'}{\varrho},$$

which yields (23). Therefore $u_h \to u$ strongly in $H_V(\Omega)$ and the proof is complete. $\square$

## 3. Energy estimates and concentration

Let us now introduce the functional $J_0 : H^1(\mathbb{R}^N) \to \mathbb{R}$ defined by

$$J_0(u) := \frac{1}{2} \sum_{i,j=1}^N \int_{\mathbb{R}^N} a_{ij}(x_0, u) D_i u D_j u + \frac{1}{2} \int_{\mathbb{R}^N} V(x_0)\, u^2 - \int_{\mathbb{R}^N} F(u)$$

where $x_0$ is as in (11). Let us set

$$\bar{c} := \inf_{\gamma \in \mathscr{P}_0} \sup_{t \in [0,1]} J_0(\gamma(t)),$$

where $\mathscr{P}_0$ is the family defined by

(27) $$\mathscr{P}_0 := \left\{ \gamma \in C([0,1], H_V(\mathbb{R}^N)) : \gamma(0) = 0,\ J_0(\gamma(1)) < 0 \right\}.$$

Let us also set

(28) $$\mathscr{P}_\varepsilon := \left\{ \gamma \in C([0,1], H_V(\Omega)) : \gamma(0) = 0,\ J_\varepsilon(\gamma(1)) < 0 \right\}.$$

In the following, if necessary, we will assume that, for every $\gamma \in \mathscr{P}_\varepsilon$, for every $t \in [0,1]$ the map $\gamma(t)$ is extended to zero outside $\Omega$.

In the next Lemma we get a critical point $u_\varepsilon$ of $J_\varepsilon$ with a precise energy upper bound.

**Lemma 3.1.** *For $\varepsilon > 0$ sufficiently small $J_\varepsilon$ admits a critical point $u_\varepsilon \in H_V(\Omega)$ such that*

(29) $$J_\varepsilon(u_\varepsilon) \leqslant \varepsilon^N \bar{c} + o(\varepsilon^N).$$

*Proof.* Let $\varepsilon > 0$. By Lemma 2.3 the functional $J_\varepsilon$ satisfies the concrete Palais–Smale condition at every level $c \in \mathbb{R}$. Moreover, since $g(x,s) = o(s)$ as $s \to 0$ uniformly in $x$, it is readily seen that there exist $\varrho_\varepsilon > 0$ and $\nu_\varepsilon > 0$ such that $J_\varepsilon$ verifies condition (63). Finally, if $z$ is a positive function in $H_V(\Omega) \setminus \{0\}$ such that $\mathrm{supt}(z) \subset \Lambda$, by (6) it results $J_\varepsilon(tz) \to -\infty$ as $t \to +\infty$. Therefore, by Proposition 6.9, minimaxing over the family (28), the functional $J_\varepsilon$ admits a nontrivial critical point $u_\varepsilon \in H_V(\Omega)$ such that

$$J_\varepsilon(u_\varepsilon) = \inf_{\gamma \in \mathscr{P}_\varepsilon} \sup_{t \in [0,1]} J_\varepsilon(\gamma(t)).$$

Since $\bar{c}$ is the Mountain–Pass value of the limiting functional $J_0$, for every $\delta > 0$ there exists a continuous path $\gamma : [0,1] \to H_V(\mathbb{R}^N)$ such that

$$(30) \qquad \bar{c} \leqslant \sup_{t \in [0,1]} J_0(\gamma(t)) \leqslant \bar{c} + \delta, \quad \gamma(0) = 0, \quad J_0(\gamma(1)) < 0.$$

Let $\zeta \in C_c^\infty(\mathbb{R}^N)$ be a cut–off function with $\zeta = 1$ in a neighbourhood $U$ of $x_0$ in $\Lambda$. We define the continuous path $\Gamma_\varepsilon : [0,1] \to H_V(\Omega)$ by setting $\Gamma_\varepsilon(\tau)(x) := \zeta(x)\gamma(\tau)\left(\frac{x-x_0}{\varepsilon}\right)$ for every $\tau \in [0,1]$ and $x \in \Omega$. Then, for every $\tau \in [0,1]$, after extension to zero outside $\Omega$, we have

$$J_\varepsilon(\Gamma_\varepsilon(\tau)) = \frac{\varepsilon^2}{2} \sum_{i,j=1}^N \int_{\mathbb{R}^N} a_{ij}\left(x, \zeta(x)\gamma(\tau)\left(\frac{x-x_0}{\varepsilon}\right)\right) D_i\zeta D_j\zeta \, \gamma^2(\tau)\left(\frac{x-x_0}{\varepsilon}\right)$$

$$+ \frac{1}{2} \sum_{i,j=1}^N \int_{\mathbb{R}^N} a_{ij}\left(x, \zeta(x)\gamma(\tau)\left(\frac{x-x_0}{\varepsilon}\right)\right) (D_i\gamma(\tau))\left(\frac{x-x_0}{\varepsilon}\right)(D_j\gamma(\tau))\left(\frac{x-x_0}{\varepsilon}\right) \zeta^2$$

$$+ \varepsilon \sum_{i,j=1}^N \int_{\mathbb{R}^N} a_{ij}\left(x, \zeta(x)\gamma(\tau)\left(\frac{x-x_0}{\varepsilon}\right)\right) D_i\zeta(D_j\gamma(\tau))\left(\frac{x-x_0}{\varepsilon}\right) \zeta\gamma(\tau)\left(\frac{x-x_0}{\varepsilon}\right)$$

$$+ \frac{1}{2} \int_{\mathbb{R}^N} V(x)\zeta^2(x)\gamma^2(\tau)\left(\frac{x-x_0}{\varepsilon}\right) - \int_{\mathbb{R}^N} G\left(x, \zeta(x)\gamma(\tau)\left(\frac{x-x_0}{\varepsilon}\right)\right).$$

Then, after the change of coordinates, for every $\tau \in [0,1]$ we get

$$J_\varepsilon(\Gamma_\varepsilon(\tau)) =$$

$$= \frac{\varepsilon^{N+2}}{2} \sum_{i,j=1}^N \int_{\mathbb{R}^N} a_{ij}\left(\varepsilon y + x_0, \zeta(\varepsilon y + x_0)\gamma(\tau)(y)\right) D_i\zeta(\varepsilon y + x_0) D_j\zeta(\varepsilon y + x_0) \gamma^2(\tau)(y)$$

$$+ \varepsilon^{N+1} \sum_{i,j=1}^N \int_{\mathbb{R}^N} a_{ij}\left(\varepsilon y + x_0, \zeta(\varepsilon y + x_0)\gamma(\tau)(y)\right) D_i\zeta(\varepsilon y + x_0) D_j\gamma(\tau)(y)\zeta(\varepsilon y + x_0)\gamma(\tau)(y)$$

$$+ \frac{\varepsilon^N}{2} \sum_{i,j=1}^N \int_{\mathbb{R}^N} a_{ij}\left(\varepsilon y + x_0, \zeta(\varepsilon y + x_0)\gamma(\tau)(y)\right) D_i\gamma(\tau)(y) D_j\gamma(\tau)(y) \zeta^2(\varepsilon y + x_0)$$

$$+ \frac{\varepsilon^N}{2} \int_{\mathbb{R}^N} V(\varepsilon y + x_0)\zeta^2(\varepsilon y + x_0)\gamma^2(\tau)(y) - \varepsilon^N \int_{\mathbb{R}^N} G(\varepsilon y + x_0, \zeta(\varepsilon y + x_0)\gamma(\tau)(y)).$$

Taking into account that for every $\tau \in [0,1]$ we have

$$\lim_{\varepsilon \to 0} \int_{\mathbb{R}^N} V(\varepsilon y + x_0)\zeta^2(\varepsilon y + x_0)\gamma^2(\tau)(y) = \int_{\mathbb{R}^N} V(x_0)\gamma^2(\tau)(y),$$

$$\lim_{\varepsilon \to 0} \int_{\mathbb{R}^N} G(\varepsilon y + x_0, \zeta(\varepsilon y + x_0)\gamma(\tau)(y)) = \int_{\mathbb{R}^N} F(\gamma(\tau)(y)),$$

and also

$$\lim_{\varepsilon \to 0} \sum_{i,j=1}^{N} \int_{\mathbb{R}^N} a_{ij}(\varepsilon y + x_0, \zeta(\varepsilon y + x_0)\gamma(\tau)(y))D_i\gamma(\tau)(y)D_j\gamma(\tau)(y)\zeta^2(\varepsilon y + x_0)$$

$$= \sum_{i,j=1}^{N} \int_{\mathbb{R}^N} a_{ij}(x_0, \gamma(\tau)(y))D_i\gamma(\tau)(y)D_j\gamma(\tau)(y),$$

we obtain

$$J_\varepsilon(\Gamma_\varepsilon(\tau)) = \varepsilon^N \left\{ \frac{1}{2} \sum_{i,j=1}^{N} \int_{\mathbb{R}^N} a_{ij}(x_0, \gamma(\tau)(y))D_i\gamma(\tau)(y)D_j\gamma(\tau)(y) \right.$$

$$\left. + \frac{1}{2} \int_{\mathbb{R}^N} V(x_0)\gamma^2(\tau)(y) - \int_{\mathbb{R}^N} F(\gamma(\tau)(y)) \right\} + o(\varepsilon^N)$$

as $\varepsilon \to 0$, namely

(31) $$J_\varepsilon(\Gamma_\varepsilon(\tau)) = \varepsilon^N J_0(\gamma(\tau)) + o(\varepsilon^N)$$

as $\varepsilon \to 0$, where $o(\varepsilon^N)$ is independent of $\tau$ (by a compactness argument). Then, by (30) and (31), it follows that $\Gamma_\varepsilon \in \mathscr{P}_\varepsilon$ for every $\varepsilon > 0$ sufficiently small and,

$$J_\varepsilon(u_\varepsilon) = \inf_{\gamma \in \mathscr{P}_\varepsilon} \sup_{t \in [0,1]} J_\varepsilon(\gamma(t)) \leqslant \sup_{t \in [0,1]} J_\varepsilon(\Gamma_\varepsilon(t))$$

$$= \varepsilon^N \sup_{t \in [0,1]} J_0(\gamma(t)) + o(\varepsilon^N)$$

$$\leqslant \varepsilon^N \bar{c} + o(\varepsilon^N) + \delta\varepsilon^N \quad \text{for every } \delta > 0.$$

By the arbitrariness of $\delta$ one concludes the proof. $\square$

In the following result we get some apriori estimates for the rescalings of $u_\varepsilon$.

**Corollary 3.2.** *Let $(\varepsilon_h) \subset \mathbb{R}^+$, $(x_h) \subset \Lambda$ and assume that $(u_{\varepsilon_h}) \subset H_V(\Omega)$ is as in Lemma 3.1. Let us set*

$$v_h \in H_V(\Omega_h), \quad \Omega_h := \varepsilon_h^{-1}(\Omega - x_h), \quad v_h(x) := u_{\varepsilon_h}(x_h + \varepsilon_h x)$$

*and put $v_h = 0$ outside $\Omega_h$.*

*Then there exists a positive constant $C$ such that*

(32) $$\|v_h\|_{H^1(\mathbb{R}^N)} \leqslant C$$

*for every $h \in \mathbb{N}$.*

*Proof.* We consider the functional $J_h : H_V(\Omega_h) \to \mathbb{R}$ given by

$$J_h(v) := \frac{1}{2} \sum_{i,j=1}^{N} \int_{\Omega_h} a_{ij}(x_h + \varepsilon_h x, v) D_i v D_j v \tag{33}$$

$$+ \frac{1}{2} \int_{\Omega_h} V(x_h + \varepsilon_h x) v^2 - \int_{\Omega_h} G(x_h + \varepsilon_h x, v).$$

Since $J_h(v_h) = \varepsilon_h^{-N} J_{\varepsilon_h}(u_{\varepsilon_h})$, by virtue of Lemma 3.1 we have $J_h(v_h) \leq \bar{c} + o(1)$ as $h \to +\infty$. Therefore, if we set $\Lambda_h = \varepsilon_h^{-1}(\Lambda - x_h)$, from inequalities (18) and (19), we get

$$\frac{\vartheta}{2} \sum_{i,j=1}^{N} \int_{\mathbb{R}^N} a_{ij}(x_h + \varepsilon_h x, v_h) D_i v_h D_j v_h + \frac{\vartheta}{2} \int_{\mathbb{R}^N} V(x_h + \varepsilon_h x) v_h^2 \tag{34}$$

$$\leq \int_{\Lambda_h} g(x_h + \varepsilon_h x, v_h) v_h + \frac{\vartheta}{2k} \int_{\mathbb{R}^N \setminus \Lambda_h} V(x_h + \varepsilon_h x) v_h^2 + \vartheta \bar{c} + o(1)$$

as $h \to +\infty$. Moreover, since by Proposition 6.4 it results $J_h'(v_h)(v_h) = 0$ for every $h \in \mathbb{N}$, again by (19), we get

$$\sum_{i,j=1}^{N} \int_{\mathbb{R}^N} a_{ij}(x_h + \varepsilon_h x, v_h) D_i v_h D_j v_h + \frac{1}{2} \sum_{i,j=1}^{N} \int_{\mathbb{R}^N} D_s a_{ij}(x_h + \varepsilon_h x, v_h) v_h D_i v_h D_j v_h$$

$$+ \int_{\mathbb{R}^N} V(x_h + \varepsilon_h x) v_h^2 \geq \int_{\Lambda_h} g(x_h + \varepsilon_h x, v_h) v_h,$$

which, in view of (9), yields

$$\left(\frac{\gamma}{2} + 1\right) \sum_{i,j=1}^{N} \int_{\mathbb{R}^N} a_{ij}(x_h + \varepsilon_h x, v_h) D_i v_h D_j v_h \tag{35}$$

$$+ \int_{\mathbb{R}^N} V(x_h + \varepsilon_h x) v_h^2 \geq \int_{\Lambda_h} g(x_h + \varepsilon_h x, v_h) v_h.$$

Then, recalling (7) and (8), by combining inequality (34) and (35) one gets

$$\min\left\{\left(\frac{\vartheta}{2} - \frac{\gamma}{2} - 1\right)\nu, \left(\frac{\vartheta}{2} - \frac{\vartheta}{2k} - 1\right)\alpha\right\} \int_{\mathbb{R}^N} \left(|Dv_h|^2 + v_h^2\right) \leq \vartheta \bar{c} + o(1) \tag{36}$$

as $h \to +\infty$, which yields the assertion. $\square$

**Corollary 3.3.** *Assume that* $(u_\varepsilon)_{\varepsilon > 0} \subset H_V(\Omega)$ *is as in Lemma 3.1.*
*Then we have*

$$\lim_{\varepsilon \to 0} \|u_\varepsilon\|_{H_V(\Omega)} = 0.$$

*Proof.* We may argue as in Step I of Lemma 2.3 with $u_h$ replaced by $u_\varepsilon$ and $c$ replaced by $J_\varepsilon(u_\varepsilon)$. Thus, from inequality (22), for every $\varepsilon > 0$ we get

$$\int_\Omega \left(|Du_\varepsilon|^2 + V(x) u_\varepsilon^2\right) \leq \frac{\vartheta}{\min\left\{\left(\frac{\vartheta}{2} - \frac{\gamma}{2} - 1\right)\nu \varepsilon^2, \frac{\vartheta}{2} - \frac{\vartheta}{2k} - 1\right\}} J_\varepsilon(u_\varepsilon).$$

By virtue of Lemma 3.1, this yields

$$\int_\Omega \left(|Du_\varepsilon|^2 + V(x)u_\varepsilon^2\right) \leqslant \frac{2\vartheta \bar{c}}{(\vartheta - \gamma - 2)\nu} \varepsilon^{N-2} + o(\varepsilon^{N-2})$$

for every $\varepsilon$ sufficiently small, which implies the assertion. $\square$

Let $\mathscr{L} : \mathbb{R}^N \times \mathbb{R} \times \mathbb{R}^N \to \mathbb{R}$ be a function of class $C^1$ such that the function $\nabla_\xi \mathscr{L}$ is of class $C^1$ and let $\varphi \in L^\infty_{\text{loc}}(\mathbb{R}^N)$. We now recall the Pucci–Serrin variational identity [PS].

**Lemma 3.4.** *Let $u : \mathbb{R}^N \to \mathbb{R}$ be a $C^2$ solution of*

$$-\text{div}\,(D_\xi \mathscr{L}(x, u, Du)) + D_s \mathscr{L}(x, u, Du) = \varphi \quad \text{in } \mathscr{D}'(\mathbb{R}^N).$$

*Then we have*

(37)
$$\sum_{i,j=1}^N \int_{\mathbb{R}^N} D_i h^j D_{\xi_i} \mathscr{L}(x, u, Du) D_j u +$$
$$- \int_{\mathbb{R}^N} \left[(\text{div}\,h)\,\mathscr{L}(x, u, Du) + h \cdot D_x \mathscr{L}(x, u, Du)\right] = \int_{\mathbb{R}^N} (h \cdot Du)\varphi$$

*for every $h \in C_c^1(\mathbb{R}^N, \mathbb{R}^N)$.*

We refer also to [DMS], where the above variational relation is proved for $C^1$ solutions.

We now derive an important consequence of the previous identity which will play an important role in the proof of Lemma 3.6.

**Lemma 3.5.** *Let $\mu > 0$ and $h, H : \mathbb{R}^+ \to \mathbb{R}$ be the continuous functions defined by*

$$h(s) = -\mu s + f(s), \qquad H(s) = \int_0^s h(t)\,dt,$$

*where $f$ satisfies (5) and (6). Moreover, let $b_{ij} \in C^1(\mathbb{R}^+) \cap L^\infty(\mathbb{R}^+)$ with $b'_{ij} \in L^\infty(\mathbb{R}^+)$ and assume that there exist $\nu' > 0$ and $R' > 0$ with*

(38)
$$\sum_{i,j=1}^N b_{ij}(s)\xi_i\xi_j \geqslant \nu'|\xi|^2, \qquad s \geqslant R' \implies \sum_{i,j=1}^N b'_{ij}(s)\xi_i\xi_j \geqslant 0$$

*for every $s \in \mathbb{R}^+$ and $\xi \in \mathbb{R}^N$.*

*Let $u \in H^1(\mathbb{R}^N)$ be any nontrivial positive solution of the equation*

(39)
$$-\sum_{j=1}^N D_j(b_{ij}(u)D_i u) + \frac{1}{2}\sum_{i,j=1}^N b'_{ij}(u)D_i u D_j u = h(u) \quad \text{in } \mathbb{R}^N.$$

*We denote by $\widehat{J}$ the associated functional*

(40)
$$\widehat{J}(u) := \frac{1}{2}\sum_{i,j=1}^N \int_{\mathbb{R}^N} b_{ij}(u)D_i u D_j u - \int_{\mathbb{R}^N} H(u).$$

*Then it results $\widehat{J}(u) \geqslant b$, where*

$$b := \inf_{\gamma \in \widehat{\mathscr{P}}} \sup_{t \in [0,1]} \widehat{J}(\gamma(t)),$$

$$\widehat{\mathscr{P}} := \left\{ \gamma \in C([0,1], H^1(\mathbb{R}^N)) : \gamma(0) = 0, \ \widehat{J}(\gamma(1)) < 0 \right\}.$$

*Proof.* By condition (38), it results

$$\widehat{J}(v) \geqslant \frac{1}{2} \min\{\nu', \mu\} \|v\|^2_{H^1(\mathbb{R}^N)} - \int_{\mathbb{R}^N} F(v) \quad \text{for every } v \in H^1(\mathbb{R}^N).$$

Then, since for every $\varepsilon > 0$ there exists $C_\varepsilon > 0$ with

$$0 \leqslant F(s) \leqslant \varepsilon s^2 + C_\varepsilon |s|^{\frac{2N}{N-2}} \quad \text{for every } s \in \mathbb{R}^+,$$

it is readily seen that there exist $\varrho_0 > 0$ and $\delta_0 > 0$ such that $\widehat{J}(v) \geqslant \delta_0$ for every $v$ with $\|v\|_{1,2} = \varrho_0$. In particular $\widehat{J}$ has a Mountain–Pass geometry. As we will see, $\widehat{\mathscr{P}} \neq \emptyset$, so that $b$ is well defined. Let $u$ be a nontrivial positive solution of (39) and consider the dilation path

$$\gamma(t)(x) := \begin{cases} u\left(\frac{x}{t}\right) & \text{if } t > 0 \\ 0 & \text{if } t = 0. \end{cases}$$

Notice that $\|\gamma(t)\|^2_{H^1} = t^{N-2}\|Du\|^2_2 + t^N \|u\|^2_2$ for every $t \in \mathbb{R}^+$, which implies that the curve $\gamma$ belongs to $C([0, +\infty[, H^1(\mathbb{R}^N))$. For every $t \in \mathbb{R}^+$ it results that

$$\widehat{J}(\gamma(t)) = \frac{1}{2} \sum_{i,j=1}^N \int_{\mathbb{R}^N} b_{ij}(\gamma(t)) D_i \gamma(t) D_j \gamma(t) - \int_{\mathbb{R}^N} H(\gamma(t))$$

$$= \frac{t^{N-2}}{2} \sum_{i,j=1}^N \int_{\mathbb{R}^N} b_{ij}(u) D_i u D_j u - t^N \int_{\mathbb{R}^N} H(u)$$

which yields, for every $t \in \mathbb{R}^+$

$$(41) \qquad \frac{d}{dt} \widehat{J}(\gamma(t)) = \frac{N-2}{2} t^{N-3} \sum_{i,j=1}^N \int_{\mathbb{R}^N} b_{ij}(u) D_i u D_j u - N t^{N-1} \int_{\mathbb{R}^N} H(u).$$

By (38), arguing like at the end of Step I of Lemma 3.6 (namely using the local Serrin estimates) it results that $u \in L^\infty_{\text{loc}}(\mathbb{R}^N)$. Hence by the regularity results of [LU], it follows that $u$ is of class $C^2$. Then we can use Lemma 3.4 by choosing $\varphi = 0$,

$$(42) \qquad \mathscr{L}(s, \xi) := \frac{1}{2} \sum_{i,j=1}^N b_{ij}(s) \xi_i \xi_j - H(s) \quad \text{for every } s \in \mathbb{R}^+ \text{ and } \xi \in \mathbb{R}^N,$$

$$h(x) := h_k(x) = T\left(\frac{x}{k}\right) x \quad \text{for every } x \in \mathbb{R}^N \text{ and } k \geqslant 1,$$

being $T \in C_c^1(\mathbb{R}^N)$ such that $T(x) = 1$ if $|x| \leqslant 1$ and $T(x) = 0$ if $|x| \geqslant 2$. In particular, it results that $h_k \in C_c^1(\mathbb{R}^N, \mathbb{R}^N)$ for every $k \geqslant 1$ and

$$D_i h_k^j(x) = D_i T\left(\frac{x}{k}\right)\frac{x_j}{k} + T\left(\frac{x}{k}\right)\delta_{ij} \quad \text{for every } x \in \mathbb{R}^N \text{ and } i,j = 1,\ldots,N$$

$$(\operatorname{div} h_k)(x) = DT\left(\frac{x}{k}\right)\cdot\frac{x}{k} + NT\left(\frac{x}{k}\right) \quad \text{for every } x \in \mathbb{R}^N.$$

Then, since $D_x \mathscr{L}(u, Du) = 0$, it follows by (37) that for every $k \geqslant 1$

$$\sum_{i,j=1}^n \int_{\mathbb{R}^N} D_i T\left(\frac{x}{k}\right)\frac{x_j}{k} D_j u D_{\xi_i}\mathscr{L}(u, Du) + \int_{\mathbb{R}^N} T\left(\frac{x}{k}\right) D_\xi \mathscr{L}(u, Du) \cdot Du$$
$$- \int_{\mathbb{R}^N} DT\left(\frac{x}{k}\right)\cdot\frac{x}{k}\mathscr{L}(u, Du) - \int_{\mathbb{R}^N} NT\left(\frac{x}{k}\right)\mathscr{L}(u, Du) = 0.$$

Since there exists $C > 0$ with

$$D_i T\left(\frac{x}{k}\right)\frac{x_j}{k} \leqslant C \quad \text{for every } x \in \mathbb{R}^N, k \geqslant 1 \text{ and } i,j = 1,\ldots,N,$$

by the Dominated Convergence Theorem, letting $k \to +\infty$, we obtain

$$\int_{\mathbb{R}^N}\left[N\mathscr{L}(u, Du) - D_\xi\mathscr{L}(u, Du)\cdot Du\right] = 0,$$

namely, by (42),

(43) $$\frac{N-2}{2}\sum_{i,j=1}^N \int_{\mathbb{R}^N} b_{ij}(u) D_i u D_j u = N\int_{\mathbb{R}^N} H(u).$$

By plugging this formula into (41), we obtain

$$\frac{d}{dt}\widehat{J}(\gamma(t)) = N(1-t^2)t^{N-3}\int_{\mathbb{R}^N} H(u)$$

which yields $\frac{d}{dt}\widehat{J}(\gamma(t)) > 0$ for $t < 1$ and $\frac{d}{dt}\widehat{J}(\gamma(t)) < 0$ for $t > 1$, i.e.

$$\sup_{t\in[0,L]}\widehat{J}(\gamma(t)) = \widehat{J}(\gamma(1)) = \widehat{J}(u).$$

Moreover, observe that

$$\gamma(0) = 0 \text{ and } \widehat{J}(\gamma(T)) < 0 \text{ for } T > 0 \text{ sufficiently large.}$$

Then, after a suitable scale change in $t$, $\gamma \in \widehat{\mathscr{P}}$ and the assertion follows. $\square$

The following is one of the main tools of the paper.

**Lemma 3.6.** *Assume that $(u_\varepsilon)_{\varepsilon>0} \subset H_V(\Omega)$ is as in Lemma 3.1.*
*Then we have*

(44) $$\lim_{\varepsilon\to 0}\max_{\partial\Lambda} u_\varepsilon = 0.$$

*Proof.* The following auxiliary fact is sufficient to prove assertion (44): if $\varepsilon_h \to 0$ and $(x_h) \subset \Lambda$ are such that $u_{\varepsilon_h}(x_h) \geqslant c$ for some $c > 0$, then

$$\lim_h V(x_h) = \min_\Lambda V. \tag{45}$$

Indeed, assume by contradiction that there exist $(\varepsilon_h) \subset \mathbb{R}^+$ with $\varepsilon_h \to 0$ and $(x_h) \subset \partial \Lambda$ such that $u_{\varepsilon_h}(x_h) \geqslant c$ for some $c > 0$. Up to a subsequence, we have $x_h \to \bar{x} \in \partial \Lambda$. Then by (45) it results

$$\min_{\partial \Lambda} V \leqslant V(\bar{x}) = \lim_h V(x_h) = \min_\Lambda V$$

which contradicts assumption (11).

We divide the proof of (45) into four steps:

<u>STEP I</u>. Up to a subsequence, $x_h \to \widehat{x}$ for some $\widehat{x} \in \Lambda$. By contradiction, we assume that

$$V(\widehat{x}) > \min_\Lambda V = V(x_0).$$

Since for every $h \in \mathbb{N}$ the function $u_{\varepsilon_h}$ solves $(P_{\varepsilon_h})$, the sequence

$$v_h \in H_V(\Omega_h), \quad \Omega_h = \varepsilon_h^{-1}(\Omega - x_h), \quad v_h(x) = u_{\varepsilon_h}(x_h + \varepsilon_h x)$$

satisfies

$$-\sum_{i,j=1}^N D_j(a_{ij}(x_h + \varepsilon_h x, v_h)D_i v_h) + \frac{1}{2}\sum_{i,j=1}^N D_s a_{ij}(x_h + \varepsilon_h x, v_h)D_i v_h D_j v_h = w_h \quad \text{in } \Omega_h,$$

$v_h > 0$ in $\Omega_h$ and $v_h = 0$ on $\partial \Omega_h$, where we have set

$$w_h := g(x_h + \varepsilon_h x, v_h) - V(x_h + \varepsilon_h x)v_h \quad \text{for every } h \in \mathbb{N}.$$

Setting $v_h = 0$ outside $\mathbb{R}^N$, by Corollary 3.2, up to a subsequence, $v_h \to v$ weakly in $H^1(\mathbb{R}^N)$. Notice that the sequence $(\chi_\Lambda(x_h + \varepsilon_h x))$ converges weak* in $L^\infty$ to a measurable function $0 \leqslant \chi \leqslant 1$. In particular, taking into account that $|w_h| \leqslant c_1|v_h| + c_2|v_h|^p$, $(w_h)$ is strongly convergent in $H^{-1}(\widetilde{\Omega})$ for every bounded subset $\widetilde{\Omega}$ of $\mathbb{R}^N$. Therefore, by a simple variant of Lemma 2.2, we conclude that $(v_h)$ is strongly convergent to $v$ in $H^1(\widetilde{\Omega})$ for every bounded subset $\widetilde{\Omega} \subset \mathbb{R}^N$ (actually, as we will see, $v_h \to v$ uniformly over compacts). Then it follows that the limit $v$ is a solution of the equation

$$-\sum_{i,j=1}^N D_j(a_{ij}(\widehat{x}, v)D_i v) + \frac{1}{2}\sum_{i,j=1}^N D_s a_{ij}(\widehat{x}, v)D_i v D_j v + V(\widehat{x})v = g_0(x, v) \quad \text{in } \mathbb{R}^N \tag{46}$$

where $g_0(x, s) := \chi(x)f(s) + (1 - \chi(x))\widetilde{f}(s)$ for every $x \in \mathbb{R}^N$ and $s \in \mathbb{R}^+$.

We now prove that $v \neq 0$. Let us set

$$d_h(x) := \begin{cases} V(x_h + \varepsilon_h x) - \frac{g(x, v_h(x))}{v_h(x)} & \text{if } v_h(x) \neq 0 \\ 0 & \text{if } v_h(x) = 0, \end{cases}$$

$$A_j(x, s, \xi) := \sum_{i=1}^N a_{ij}(x_h + \varepsilon_h x, s)\xi_i \quad \text{for } j = 1, \ldots, N,$$

$$B(x, s, \xi) := d_h(x)s,$$

$$C(x, s) := \frac{1}{2} \sum_{i,j=1}^N D_s a_{ij}(x_h + \varepsilon_h x, s) D_i v_h(x) D_j v_h(x)$$

for every $x \in \mathbb{R}^N$, $s \in \mathbb{R}^+$ and $\xi \in \mathbb{R}^N$. Taking into account the assumptions on the coefficients $a_{ij}(x, s)$, it results that

$$A(x, s, \xi) \cdot \xi \geqslant \nu |\xi|^2, \quad |A(x, s, \xi)| \leqslant c|\xi|, \quad |B(x, s, \xi)| \leqslant d_h(x)|s|.$$

Moreover, by (10) we have

$$s \geqslant R \implies C(x, s)s \geqslant 0$$

for every $x \in \mathbb{R}^N$ and $s \in \mathbb{R}^+$. By the growth condition on $g$, $d_h \in L^{\frac{N}{2-\delta}}(B_{2\varrho}(0))$ for every $\varrho > 0$ and

$$S = \sup_h \|d_h\|_{L^{\frac{N}{2-\delta}}(B_{2\varrho}(0))} \leqslant D_\varrho \sup_{h \in \mathbb{N}} \|v_h\|_{L^{2^*}(B_{2\varrho}(0))} < +\infty$$

for some $\delta > 0$ sufficiently small. Since we have $\text{div}(A(x, v_h, Dv_h)) = B(x, v_h, Dv_h) + C(x, v_h)$ for every $h \in \mathbb{N}$, by virtue of [SE, Theorem 1 and Remark at p.261] there exists a positive constant $M(\delta, N, c, \varrho^\delta S)$ and a radius $\varrho > 0$, sufficiently small, such that

$$\sup_{h \in \mathbb{N}} \max_{x \in B_\varrho(0)} |v_h(x)| \leqslant M(\delta, N, c, \varrho^\delta S)(2\varrho)^{-N/2} \sup_{h \in \mathbb{N}} \|v_h\|_{L^2(B_{2\varrho}(0))} < +\infty$$

so that $(v_h)$ is uniformly bounded in $B_\varrho(0)$. Then, by [SE, Theorem 8], $(v_h)$ is bounded in some $C^{1,\alpha}(\overline{B_{\varrho/2}(0)})$. Up to a subsequence this implies that $(v_h)$ converges uniformly to $v$ in $\overline{B_{\varrho/2}(0)}$. This yields $v(0) = \lim_h v_h(0) = \lim_h u_{\varepsilon_h}(x_h) \geqslant c > 0$.

In a similar fashion one shows that $v_h \to v$ uniformly over compacts.

<u>STEP II</u>. We prove that $v$ actually solves the following equation

$$(47) \quad -\sum_{i,j=1}^N D_j(a_{ij}(\widehat{x}, v)D_i v) + \frac{1}{2}\sum_{i,j=1}^N D_s a_{ij}(\widehat{x}, v)D_i v D_j v + V(\widehat{x})v = f(v) \quad \text{in } \mathbb{R}^N.$$

In general the function $\chi$ of Step I is given by $\chi = \chi_{T_\Lambda(\widehat{x})}$, where $T_\Lambda(\widehat{x})$ is the tangent cone of $\Lambda$ at $\widehat{x}$. On the other hand, since we may assume without loss of generality that $\Lambda$ is smooth, it results (up to a rotation) that $\chi(x) = \chi_{\{x_1 < 0\}}(x)$ for every $x \in \mathbb{R}^N$. In particular, $v$ is a solution of the following problem

$$
(48) \quad -\sum_{i,j=1}^{N} D_j(a_{ij}(\widehat{x},v)D_i v) + \frac{1}{2}\sum_{i,j=1}^{N} D_s a_{ij}(\widehat{x},v)D_i v D_j v \\
+ V(\widehat{x})v = \chi_{\{x_1<0\}}(x)f(v) + \chi_{\{x_1>0\}}(x)\widetilde{f}(v) \quad \text{in } \mathbb{R}^N.
$$

Let us first show that $v(x) \leqslant \ell$ on $\{x_1 = 0\}$. To this aim, let us use again Lemma 3.4, by choosing this time

$$\varphi(x) := \chi_{\{x_1<0\}}(x)f(v(x)) + \chi_{\{x_1>0\}}(x)\widetilde{f}(v(x)) \quad \text{for every } x \in \mathbb{R}^N$$

$$\mathscr{L}(s,\xi) := \frac{1}{2}\sum_{i,j=1}^{N} a_{ij}(\widehat{x},s)\xi_i\xi_j + \frac{V(\widehat{x})}{2}s^2 \quad \text{for every } s \in \mathbb{R} \text{ and } \xi \in \mathbb{R}^N,$$

$$h(x) := h_k(x) = \left(T\left(\frac{x}{k}\right), 0, \ldots, 0\right) \quad \text{for every } x \in \mathbb{R}^N \text{ and } k \geqslant 1.$$

Then $h_k \in C_c^1(\mathbb{R}^N, \mathbb{R}^N)$ and, since $D_x\mathscr{L}(v, Dv) = 0$, for every $k \geqslant 1$ it results

$$\int_{\mathbb{R}^N} \left[\frac{1}{k}\sum_{i=1}^{N} D_i T\left(\frac{x}{k}\right) D_1 v D_{\xi_i}\mathscr{L}(v, Dv) - D_1 T\left(\frac{x}{k}\right)\frac{1}{k}\mathscr{L}(v, Dv)\right]$$
$$= \int_{\mathbb{R}^N} T\left(\frac{x}{k}\right)\varphi(x, v)D_1 v.$$

Again by the Dominated Convergence Theorem, letting $k \to +\infty$, it results

$$\int_{\mathbb{R}^N} \varphi(x, v) D_1 v = 0,$$

that is, after integration by parts,

$$\int_{\mathbb{R}^{N-1}} \left[F(v(0, x')) - \widetilde{F}(v(0, x'))\right] dx' = 0.$$

Taking into account that $F(s) \geqslant \widetilde{F}(s)$ with equality only if $s \leqslant \ell$, we get $v(0, x') \leqslant \ell$ for every $x' \in \mathbb{R}^{N-1}$. To prove that actually $v(x_1, x') \leqslant \ell$ for every $x_1 > 0$ and $x' \in \mathbb{R}^{N-1}$, we test (48) with

$$\eta(x) := \begin{cases} 0 & \text{if } x_1 < 0 \\ (v(x_1, x') - \ell)^+ \exp\{\zeta(v(x_1, x'))\} & \text{if } x_1 > 0 \end{cases}$$

where $\zeta(s)$ is as in (25) and then we argue as in Section 4 (see the computations in formula (55)). In particular,

$$(49) \quad \varphi(x, v(x)) = f(v(x)) \quad \text{for every } x \in \mathbb{R}^N,$$

so that $v$ is a nontrivial solution of (47).

<u>STEP III</u> If $J_h : H_V(\Omega_h) \to \mathbb{R}$ is as in (33), the function $v_h$ is a critical point of $J_h$ and $J_h(v_h) = \varepsilon_h^{-N} J_{\varepsilon_h}(u_{\varepsilon_h})$ for every $h \in \mathbb{N}$. Let us consider the functional $J_{\widehat{x}} : H^1(\mathbb{R}^N) \to \mathbb{R}$ defined as

$$J_{\widehat{x}}(u) := \frac{1}{2} \sum_{i,j=1}^{N} \int_{\mathbb{R}^N} a_{ij}(\widehat{x}, u) D_i u D_j u + \frac{1}{2} \int_{\mathbb{R}^N} V(\widehat{x}) u^2 - \int_{\mathbb{R}^N} F(u).$$

We now want to prove that

(50) $$J_{\widehat{x}}(v) \leqslant \liminf_{h} J_h(v_h).$$

Let us set for every $h \in \mathbb{N}$ and $x \in \Omega_h$

(51) $$\xi_h(x) := \frac{1}{2} \sum_{i,j=1}^{N} a_{ij}(x_h + \varepsilon_h x, v_h) D_i v_h D_j v_h + \frac{1}{2} V(x_h + \varepsilon_h x) v_h^2 - G(x_h + \varepsilon_h x, v_h).$$

Since $v_h \to v$ in $H^1$ over compact sets, in view of (49), for every $\varrho > 0$ one gets

$$\lim_{h} \int_{B_\varrho(0)} \xi_h(x) = \frac{1}{2} \int_{B_\varrho(0)} \left( \sum_{i,j=1}^{N} a_{ij}(\widehat{x}, v) D_i v D_j v + V(\widehat{x}) v^2 \right) - \int_{B_\varrho(0)} F(v).$$

Moreover, as $v$ belongs to $H^1(\mathbb{R}^N)$,

$$\frac{1}{2} \int_{B_\varrho(0)} \left( \sum_{i,j=1}^{N} a_{ij}(\widehat{x}, v) D_i v D_j v + V(\widehat{x}) v^2 \right) - \int_{B_\varrho(0)} F(v) = J_{\widehat{x}}(v) - o(1)$$

as $\varrho \to +\infty$. Therefore, it suffices to show that for every $\delta > 0$ there exists $\varrho > 0$ with

(52) $$\liminf_{h} \int_{\Omega_h \setminus B_\varrho(0)} \xi_h(x) \geqslant -\delta.$$

Consider a function $\eta_\varrho \in C^\infty(\mathbb{R}^N)$ such that $0 \leqslant \eta_\varrho \leqslant 1$, $\eta_\varrho = 0$ on $B_{\varrho-1}(0)$, $\eta_\varrho = 1$ on $\mathbb{R}^N \setminus B_\varrho(0)$ and $|D\eta_\varrho| \leqslant c$. Let us set for every $h \in \mathbb{N}$

$$\beta_h(\varrho) := \sum_{i,j=1}^{N} \int_{B_\varrho(0) \setminus B_{\varrho-1}(0)} a_{ij}(x_h + \varepsilon_h x, v_h) D_i v_h D_j(\eta_\varrho v_h)$$

$$+ \frac{1}{2} \sum_{i,j=1}^{N} \int_{B_\varrho(0) \setminus B_{\varrho-1}(0)} D_s a_{ij}(x_h + \varepsilon_h x, v_h) \eta_\varrho v_h D_i v_h D_j v_h$$

$$+ \int_{B_\varrho(0) \setminus B_{\varrho-1}(0)} V(x_h + \varepsilon_h x) v_h^2 \eta_\varrho - \int_{B_\varrho(0) \setminus B_{\varrho-1}(0)} g(x_h + \varepsilon_h x, v_h) \eta_\varrho v_h.$$

After some computations, in view of (9), (51) and Proposition 6.4, one gets

$$-\beta_h(\varrho) + J_h'(v_h)(\eta_\varrho v_h)$$
$$\leqslant (\gamma + 2) \int_{\Omega_h \setminus B_\varrho(0)} \xi_h(x) - \frac{\gamma}{2} \int_{\Omega_h \setminus B_\varrho(0)} V(x_h + \varepsilon_h x) v_h^2$$
$$+ (\gamma + 2) \int_{\Omega_h \setminus B_\varrho(0)} G(x_h + \varepsilon_h x, v_h) - \int_{\Omega_h \setminus B_\varrho(0)} g(x_h + \varepsilon_h x, v_h) v_h$$

Notice that, by virtue of (18), for $\varrho$ large enough, setting $\Lambda_h = \varepsilon_h^{-1}(\Lambda - x_h)$, we get

$$-\frac{\gamma}{2} \int_{\Lambda_h \setminus B_\varrho(0)} V(x_h + \varepsilon_h x) v_h^2 + (\gamma + 2) \int_{\Lambda_h \setminus B_\varrho(0)} G(x_h + \varepsilon_h x, v_h)$$
$$- \int_{\Lambda_h \setminus B_\varrho(0)} g(x_h + \varepsilon_h x, v_h) v_h \leqslant -(\vartheta - 2 - \gamma) \int_{\Lambda_h \setminus B_\varrho(0)} G(x_h + \varepsilon_h x, v_h) \leqslant 0.$$

Analogously, in view of (19), we obtain

$$-\frac{\gamma}{2} \int_{\Omega_h \setminus (B_\varrho(0) \cup \Lambda_h)} V(x_h + \varepsilon_h x) v_h^2 + (\gamma + 2) \int_{\Omega_h \setminus (B_\varrho(0) \cup \Lambda_h)} G(x_h + \varepsilon_h x, v_h)$$
$$- \int_{\Omega_h \setminus (B_\varrho(0) \cup \Lambda_h)} g(x_h + \varepsilon_h x, v_h) v_h$$
$$\leqslant -\frac{\gamma}{2} \int_{\Omega_h \setminus (B_\varrho(0) \cup \Lambda_h)} V(x_h + \varepsilon_h x) v_h^2 + \frac{\gamma}{2k} \int_{\Omega_h \setminus (B_\varrho(0) \cup \Lambda_h)} V(x_h + \varepsilon_h x) v_h^2 \leqslant 0.$$

Therefore, since $J_h'(v_h)(\eta_\varrho v_h) = 0$ for every $h \in \mathbb{N}$ and

$$\limsup_h \beta_h(\varrho) = o(1) \quad \text{as } \varrho \to +\infty,$$

inequality (52) follows and thus (50) holds true.

STEP IV. In this step we get the desired contradiction. By combining Lemma 3.1 with the inequality (50), one immediately gets

(53) $$J_{\widehat{x}}(v) \leqslant \bar{c} = \inf_{\gamma \in \mathscr{P}_0} \sup_{t \in [0,1]} J_0(\gamma(t)).$$

Since $v$ is a nontrivial solution of (47), by applying Lemma 3.5 with

$$\mu = V(\widehat{x}), \quad \nu' = \nu, \quad R' = R, \quad b_{ij}(s) = a_{ij}(\widehat{x}, s),$$

being $\widehat{\mathscr{P}} \subseteq \mathscr{P}_0$, $V(\widehat{x}) > V(x_0)$ and, by (12),

$$\sum_{i,j=1}^N a_{ij}(\widehat{x}, s) \xi_i \xi_j \geqslant \sum_{i,j=1}^N a_{ij}(x_0, s) \xi_i \xi_j \quad \text{for every } s \in \mathbb{R}^+ \text{ and } \xi \in \mathbb{R}^N,$$

it follows that

(54) $$J_{\widehat{x}}(v) \geqslant \inf_{\gamma \in \widehat{\mathscr{P}}} \sup_{t \in [0,1]} J_{\widehat{x}}(\gamma(t)) > \inf_{\gamma \in \mathscr{P}_0} \sup_{t \in [0,1]} J_0(\gamma(t)) = \bar{c},$$

which contradicts (53). □

## 4. Proof of the main result

We are now ready to prove Theorem 1.1.

<u>Step I</u>. We prove that $(a)$ holds. By Lemma 3.6 there exists $\varepsilon_0 > 0$ such that
$$u_\varepsilon(x) < \ell \quad \text{for every } \varepsilon \in (0, \varepsilon_0) \text{ and } x \in \partial \Lambda.$$
Then, since $u_\varepsilon \in H_V(\Omega)$, for every $\varepsilon \in (0, \varepsilon_0)$, if $\zeta$ is defined as in (25), the function
$$v_\varepsilon(x) := \begin{cases} 0 & \text{if } x \in \Lambda \\ (u_\varepsilon(x) - \ell)^+ \exp\{\zeta(u_\varepsilon(x))\} & \text{if } x \in \Omega \setminus \Lambda \end{cases}$$
belongs to $H_0^1(\Omega)$ and (by Proposition 6.4) it is an admissible test for the equation
$$-\varepsilon^2 \sum_{i,j=1}^N D_j(a_{ij}(x, u_\varepsilon) D_i u_\varepsilon) + \frac{\varepsilon^2}{2} \sum_{i,j=1}^N D_s a_{ij}(x, u_\varepsilon) D_i u_\varepsilon D_j u_\varepsilon + V(x) u_\varepsilon = g(x, u_\varepsilon).$$
After some computations, one obtains
$$(55) \quad \varepsilon^2 \sum_{i,j=1}^N \int_{\Omega \setminus \Lambda} a_{ij}(x, u_\varepsilon) D_i\big[(u_\varepsilon - \ell)^+\big] D_j\big[(u_\varepsilon - \ell)^+\big] \exp\{\zeta(u_\varepsilon)\}$$
$$+ \varepsilon^2 \sum_{i,j=1}^N \int_{\Omega \setminus \Lambda} \left[\frac{1}{2} D_s a_{ij}(x, u_\varepsilon) + \zeta'(u_\varepsilon) a_{ij}(x, u_\varepsilon)\right] D_i u_\varepsilon D_j u_\varepsilon (u_\varepsilon - \ell)^+ \exp\{\zeta(u_\varepsilon)\}$$
$$+ \int_{\Omega \setminus \Lambda} \Phi_\varepsilon(x) \big[(u_\varepsilon - \ell)^+\big]^2 \exp\{\zeta(u_\varepsilon)\} + \int_{\Omega \setminus \Lambda} \Phi_\varepsilon(x) \ell (u_\varepsilon - \ell)^+ \exp\{\zeta(u_\varepsilon)\} = 0,$$
where $\Phi_\varepsilon : \Omega \to \mathbb{R}$ is the function given by
$$\Phi_\varepsilon(x) := V(x) - \frac{g(x, u_\varepsilon(x))}{u_\varepsilon(x)}.$$
Notice that, by virtue of condition (19), one has
$$\Phi_\varepsilon(x) > 0 \quad \text{for every } x \in \Omega \setminus \Lambda.$$
Therefore, taking into account (26), all the terms in (55) must be equal to zero. We conclude that $(u_\varepsilon - \ell)^+ = 0$ on $\Omega \setminus \Lambda$, namely,

(56) $\qquad u_\varepsilon(x) \leqslant \ell \quad \text{for every } \varepsilon \in (0, \varepsilon_0) \text{ and } x \in \Omega \setminus \Lambda.$

Hence, by Proposition 2.1, $u_\varepsilon$ is a positive solution of the original problem $(P_\varepsilon)$. Moreover, by virtue of (10), using again the argument at the end of Step I of Lemma 3.6 it results that $u_\varepsilon \in L^\infty_{\text{loc}}(\Omega)$, which, by the regularity results of [LU], yields $u_\varepsilon \in C(\overline{\Omega})$. Notice that by arguing in a similar fashion testing with
$$v_\varepsilon(x) := \begin{cases} 0 & \text{if } x \in \Lambda \\ (u_\varepsilon(x) - \sup_{\partial \Lambda} u_\varepsilon)^+ \exp\{\zeta(u_\varepsilon(x))\} & \text{if } x \in \Omega \setminus \Lambda \end{cases}$$

it results $u_\varepsilon \to 0$ uniformly outside $\Lambda$.

STEP II. We prove that (b) holds. If $x_\varepsilon$ denotes the maximum of $u_\varepsilon$ in $\Lambda$, since $u_\varepsilon \to 0$ uniformly outside $\Lambda$, it results that $u_\varepsilon(x_\varepsilon) = \sup_\Omega u_\varepsilon$. By arguing as at the end of Step I of Lemma 3.6, setting $v_\varepsilon(x) = u_\varepsilon(x_\varepsilon + \varepsilon x)$ it results that the sequence $(v_\varepsilon(0))$ is bounded in $\mathbb{R}$. Then there exists $\sigma' > 0$ such that $u_\varepsilon(x_\varepsilon) = v_\varepsilon(0) \leqslant \sigma'$. Assume now by contradiction that $u_\varepsilon(x_\varepsilon) \leqslant \sigma$ for some $\varepsilon \in (0, \varepsilon_0)$. Then, taking into account the definition of $\sigma$ and that $u_\varepsilon \to 0$ uniformly outside $\Lambda$, it holds (with strict inequality in some subset of $\Omega$)

$$(57) \qquad V(x) - \frac{f(u_\varepsilon(x))}{u_\varepsilon(x)} \geqslant 0 \quad \text{for every } x \in \Omega.$$

Let $\zeta : \mathbb{R}^+ \to \mathbb{R}$ be the map defined in (25). Then, in view of Proposition 6.4, the function $u_\varepsilon \exp\{\zeta(u_\varepsilon)\}$ can be chosen as an admissible test in the equation

$$-\varepsilon^2 \sum_{i,j=1}^N D_j(a_{ij}(x, u_\varepsilon)D_i u_\varepsilon) + \frac{\varepsilon^2}{2} \sum_{i,j=1}^N D_s a_{ij}(x, u_\varepsilon)D_i u_\varepsilon D_j u_\varepsilon + V(x)u_\varepsilon = f(u_\varepsilon).$$

After some computations, one obtains

$$(58) \qquad \varepsilon^2 \sum_{i,j=1}^N \int_\Omega a_{ij}(x, u_\varepsilon) D_i u_\varepsilon D_j u_\varepsilon \exp\{\zeta(u_\varepsilon)\}$$

$$+ \varepsilon^2 \sum_{i,j=1}^N \int_\Omega \left[\frac{1}{2} D_s a_{ij}(x, u_\varepsilon) + \zeta'(u_\varepsilon)a_{ij}(x, u_\varepsilon)\right] D_i u_\varepsilon D_j u_\varepsilon u_\varepsilon \exp\{\zeta(u_\varepsilon)\}$$

$$+ \int_\Omega \left(V(x) - \frac{f(u_\varepsilon)}{u_\varepsilon}\right) u_\varepsilon^2 \exp\{\zeta(u_\varepsilon)\} = 0.$$

Then, by (8), (26) and (57) all the terms in equation (58) must be equal to zero, namely $u_\varepsilon \equiv 0$, which is not possible. Then $u_\varepsilon(x_\varepsilon) \geqslant \sigma$ for every $\varepsilon \in (0, \varepsilon_0)$ and by (45) we also get $d(x_\varepsilon, \mathscr{M}) \to 0$ as $\varepsilon \to 0$.

STEP III. We prove that (c) holds. Assume by contradiction that there exists $\varrho > 0$, $\delta > 0$, $\varepsilon_h \to 0$ and $y_h \in \Lambda \setminus B_\varrho(x_{\varepsilon_h})$ such that

$$(59) \qquad \limsup_h u_{\varepsilon_h}(y_h) \geqslant \delta.$$

Then, arguing as in Lemma 3.6, we can assume that $y_h \to y$, $x_{\varepsilon_h} \to \widetilde{y}$ and $v_h(y) := u_{\varepsilon_h}(y_h + \varepsilon_h y) \to v$, $\widetilde{v}_h(y) := u_{\varepsilon_h}(x_{\varepsilon_h} + \varepsilon_h y) \to \widetilde{v}$ strongly in $H^1_{\text{loc}}(\mathbb{R}^N)$, where $v$ is a solution of

$$-\sum_{i,j=1}^N D_j(a_{ij}(y,v)D_i v) + \frac{1}{2} \sum_{i,j=1}^N D_s a_{ij}(y,v) D_i v D_j v + V(y)v = f(v) \quad \text{in } \mathbb{R}^N$$

and $\widetilde{v}$ is a solution of

$$-\sum_{i,j=1}^N D_j(a_{ij}(\widetilde{y},v)D_i v) + \frac{1}{2} \sum_{i,j=1}^N D_s a_{ij}(\widetilde{y},v) D_i v D_j v + V(\widetilde{y})v = f(v) \quad \text{in } \mathbb{R}^N.$$

Observe that $v \neq 0$ and $\widetilde{v} \neq 0$. Indeed, arguing as in Step I of Lemma 3.6 it results that $(v_h)$ and $(\widetilde{v}_h)$ converge uniformly in a neighbourhood of zero, so that from (59) and $u_{\varepsilon_h}(x_{\varepsilon_h}) \geqslant \sigma$ we get $v(0) \geqslant \delta$ and $\widetilde{v}(0) \geqslant \sigma$. Now, setting $z_h := \frac{x_{\varepsilon_h} - y_h}{\varepsilon_h}$ and

$$\xi_h(y) := \frac{1}{2} \sum_{i,j=1}^{N} a_{ij}(y_h + \varepsilon_h y, v_h) D_i v_h D_j v_h + \frac{1}{2} V(y_h + \varepsilon_h y) v_h^2 - G(y_h + \varepsilon_h y, v_h),$$

if $\psi \in C^\infty(\mathbb{R})$, $0 \leqslant \psi \leqslant 1$, $\psi(s) = 0$ for $s \leqslant 1$ and $\psi(s) = 1$ for $s \geqslant 2$, arguing as in Lemma 3.6 by testing the equation satisfied by $v_h$ with

$$\varphi_{h,R}(y) = v_h(y) \left[ \psi\left(\frac{|y|}{R}\right) + \psi\left(\frac{|y - z_h|}{R}\right) - 1 \right],$$

taking into account that

$$\lim_h \left| \int_{B_{2R}(0) \cup B_{2R}(z_h) \setminus (B_R(0) \cup B_R(z_h))} \xi_h(y) \right| = o(1)$$

as $R \to +\infty$, it turns out that for every $\delta > 0$ there exists $R > 0$ with

$$\liminf_h \int_{\Omega_h \setminus (B_R(0) \cup B_R(z_h))} \xi_h(y) \geqslant -\delta.$$

Moreover, for every $R > 0$ we have

$$\liminf_h \int_{B_R(0) \cup B_R(z_h)} \xi_h(y)$$

$$= \liminf_h \int_{B_R(0)} \frac{1}{2} \sum_{i,j=1}^{N} a_{ij}(y_h + \varepsilon_h y, v_h) D_i v_h D_j v_h + \frac{1}{2} V(y_h + \varepsilon_h y) v_h^2 - G(y_h + \varepsilon_h y, v_h)$$

$$+ \liminf_h \int_{B_R(0)} \frac{1}{2} \sum_{i,j=1}^{N} a_{ij}(x_{\varepsilon_h} + \varepsilon_h y, \widetilde{v}_h) D_i \widetilde{v}_h D_j \widetilde{v}_h + \frac{1}{2} V(x_{\varepsilon_h} + \varepsilon_h y) \widetilde{v}_h^2 - G(x_{\varepsilon_h} + \varepsilon_h y, \widetilde{v}_h)$$

$$= \int_{B_R(0)} \frac{1}{2} \sum_{i,j=1}^{N} a_{ij}(y, v) D_i v D_j v + \frac{1}{2} V(y) v^2 - F(v)$$

$$+ \int_{B_R(0)} \frac{1}{2} \sum_{i,j=1}^{N} a_{ij}(\widetilde{y}, \widetilde{v}) D_i \widetilde{v} D_j \widetilde{v} + \frac{1}{2} V(\widetilde{y}) \widetilde{v}^2 - F(\widetilde{v}).$$

Therefore, we deduce that

$$\liminf_h \varepsilon_h^{-N} J_{\varepsilon_h}(u_{\varepsilon_h}) = \liminf_h \int_{\Omega_h} \xi_h(y) \geqslant J_y(v) + J_{\widetilde{y}}(\widetilde{v}).$$

If $b_y$ and $b_{\widetilde{y}}$ denote the Mountain–Pass values of the functionals $J_y$ and $J_{\widetilde{y}}$ respectively, by Lemma 3.5, (11) and (12) we have $J_y(v) \geqslant b_y \geqslant \bar{c}$ and $J_{\widetilde{y}}(\widetilde{v}) \geqslant b_{\widetilde{y}} \geqslant \bar{c}$. Therefore we conclude that $\liminf_h \varepsilon_h^{-N} J_{\varepsilon_h}(u_{\varepsilon_h}) \geqslant 2\bar{c}$, which contradicts Lemma 3.1.

STEP IV. We prove that (d) holds. By Corollary 3.3, we have $\|u_\varepsilon\|_{H_V(\Omega)} \to 0$. In particular, $u_\varepsilon \to 0$ in $L^q(\Omega)$ for every $2 \leqslant q \leqslant 2^*$. As a consequence $u_\varepsilon \to 0$ in $L^q(\Omega)$ also for every $q > 2^*$. Indeed, if $q > 2^*$, we have

$$\int_\Omega |u_\varepsilon|^q = \int_\Omega |u_\varepsilon|^{q-2^*} |u_\varepsilon|^{2^*} \leqslant \sigma'^{q-2^*} \int_\Omega |u_\varepsilon|^{2^*} \to 0$$

as $\varepsilon \to 0$.

The proof is now complete. □

## 5. A few related open problems

We quote here a few (open) problems related to the main result of the paper.

**Problem 5.1.** Under suitable assumptions, does a Gidas–Ni–Nirenberg [GNN] type result (radial symmetry) hold for the solutions of autonomous equations of the type

$$(60) \qquad -\sum_{i,j=1}^N D_j(b_{ij}(u)D_i u) + \frac{1}{2}\sum_{i,j=1}^N b'_{ij}(u)D_i u D_j u = h(u) \quad \text{in } \mathbb{R}^N?$$

**Problem 5.2.** Under suitable assumptions on $b_{ij}$ and $h$, is it possible to prove, as in the semilinear case, a *uniqueness* result for the solutions of equation (60)?

**Problem 5.3.** Under suitable assumptions on $b_{ij}$ and $h$, is it possible to prove, as in the semilinear case, that there exists a *least energy* solution of equation (60)? In other words, is there a positive solution $\omega \in H^1(\mathbb{R}^N)$ such that

$$\widetilde{J}(\omega) = \inf\left\{\widetilde{J}(u): \ u \in H^1(\mathbb{R}^N) \setminus \{0\} \text{ is a solution of } (60)\right\},$$

being $\widetilde{J}$ the functional associated with (60)? We believe so, and in particular that this solution correspond exactly to the Mountain–Pass solution.

**Problem 5.4.** Is it true that for each $\varepsilon > 0$ the solution $u_\varepsilon$ of problem $(P_\varepsilon)$ admits a *unique* maximum point inside $\Lambda$?

**Problem 5.5.** Is it true that the solutions $u_\varepsilon$ of problem $(P_\varepsilon)$ *decay exponentially* as for the semilinear case (see formula (16))?

## 6. Appendix: recalls of nonsmooth critical point theory

In this section we quote from [CDM, DM] some tools of nonsmooth critical point theory that are used in the paper.

For the sake of completeness, let us recall the definition of weak slope.

**Definition 6.1.** Let $X$ be a complete metric space, $f : X \to \mathbb{R}$ be a continuous function, and $u \in X$. We denote by $|df|(u)$ the supremum of the real numbers $\sigma \geqslant 0$ such that there exist $\delta > 0$ and a continuous map $\mathscr{H} : B(u,\delta) \times [0,\delta] \to X$ such that, for every $v$ in $B(u,\delta)$, and for every $t$ in $[0,\delta]$ it results

$$d(\mathscr{H}(v,t),v) \leqslant t, \quad f(\mathscr{H}(v,t)) \leqslant f(v) - \sigma t.$$

The extended real number $|df|(u)$ is called the weak slope of $f$ at $u$.

**Definition 6.2.** We say that $u \in X$ is a critical point of $f$ if $|df|(u) = 0$. We say that $c \in \mathbb{R}$ is a critical value of $f$ if there exists a critical point $u \in X$ of $f$ with $f(u) = c$.

**Definition 6.3.** Let $c \in \mathbb{R}$. We say that $f$ satisfies the Palais–Smale condition at level $c$ ($(PS)_c$ in short), if every sequence $(u_h)$ in $X$ such that $|df|(u_h) \to 0$ and $f(u_h) \to c$ admits a subsequence converging in $X$.

Let us now return to the concrete setting and choose $X = H_V(\Omega)$. Let $\varepsilon > 0$ and consider the functional $f : H_V(\Omega) \to \mathbb{R}$ defined by setting

$$(61) \qquad f(u) = \frac{\varepsilon^2}{2} \sum_{i,j=1}^{N} \int_\Omega a_{ij}(x,u) D_i u D_j u + \frac{1}{2} \int_\Omega V(x) u^2 - \int_\Omega G(x,u)$$

where $g : \Omega \times \mathbb{R} \to \mathbb{R}$ is now any Carathéodory map and $G(x,s) = \int_0^s g(x,t)\,dt$. Although $f$ is merely continuous, its directional derivatives exist along some special directions.

**Proposition 6.4.** *Let $u, \varphi \in H_V(\Omega)$ be such that*

$$\left[ \left( \sum_{i,j=1}^{N} D_s a_{ij}(x,u) D_i u D_j u \right) \varphi \right]^{-} \in L^1(\Omega).$$

*Then*

$$\left( \sum_{i,j=1}^{N} D_s a_{ij}(x,u) D_i u D_j u \right) \varphi \in L^1(\Omega),$$

*the directional derivative $f'(u)(\varphi)$ exists, and it holds*

$$f'(u)(\varphi) = \varepsilon^2 \sum_{i,j=1}^{N} \int_\Omega a_{ij}(x,u) D_i u D_j \varphi + \frac{\varepsilon^2}{2} \sum_{i,j=1}^{N} \int_\Omega D_s a_{ij}(x,u) D_i u D_j u \varphi$$

$$+ \int_\Omega V(x) u \varphi - \int_\Omega g(x,u) \varphi.$$

*In particular, if (10) holds, for every $\varphi \in L^\infty(\Omega)$, $\varphi \geqslant 0$ the derivative $f'(u)(\varphi u)$ exists.*

**Definition 6.5.** We say that $u$ is a weak solution of the problem

$$(62) \quad \begin{cases} -\varepsilon^2 \sum\limits_{i,j=1}^{N} D_j(a_{ij}(x,u) D_i u) + \frac{\varepsilon^2}{2} \sum\limits_{i,j=1}^{N} D_s a_{ij}(x,u) D_i u D_j u + V(x) u = g(x,u) & \text{in } \Omega \\ u = 0 & \text{on } \partial\Omega \end{cases}$$

if $u \in H_V(\Omega)$ and

$$-\varepsilon^2 \sum_{i,j=1}^{N} D_j(a_{ij}(x,u) D_i u) + \frac{\varepsilon^2}{2} \sum_{i,j=1}^{N} D_s a_{ij}(x,u) D_i u D_j u + V(x) u = g(x,u)$$

is satisfied in $\mathscr{D}'(\Omega)$.

We now introduce a variant of the classical Palais–Smale condition, suitable for our purposes.

**Definition 6.6.** Let $c \in \mathbb{R}$. We say that $(u_h) \subset H_V(\Omega)$ is a concrete Palais–Smale sequence at level $c$ $((CPS)_c$–sequence, in short) for the $f$, if $f(u_h) \to c$ and

$$\sum_{i,j=1}^N D_s a_{ij}(x, u_h) D_i u_h D_j u_h \in (H_V(\Omega))' \quad \text{as } h \to +\infty,$$

$$-\varepsilon^2 \sum_{i,j=1}^N D_j(a_{ij}(x, u_h) D_i u_h) + \frac{\varepsilon^2}{2} \sum_{i,j=1}^N D_s a_{ij}(x, u_h) D_i u_h D_j u_h + V(x) u_h - g(x, u_h) \to 0$$

strongly in $(H_V(\Omega))'$. We say that $f$ satisfies the concrete Palais–Smale condition at level $c$ $((CPS)_c$ condition), if every $(CPS)_c$–sequence for $f$ admits a strongly convergent subsequence in $H_V(\Omega)$.

The next result allows us to connect the critical points of $f$ (as in Definition 6.2) with the weak solutions of problem (62).

**Proposition 6.7.** *For every $u$ in $H_V(\Omega)$ we have*

$$|df|(u) \geqslant \sup\{\langle w_u, \varphi \rangle : \varphi \in C_c^\infty(\Omega), \|\varphi\|_{H_V} \leqslant 1\}$$

*where*

$$w_u = -\varepsilon^2 \sum_{i,j=1}^N D_j(a_{ij}(x, u) D_i u) + \frac{\varepsilon^2}{2} \sum_{i,j=1}^N D_s a_{ij}(x, u) D_i u D_j u + V(x) u - g(x, u).$$

*In particular, if $|df|(u) < +\infty$, it follows that*

$$-\varepsilon^2 \sum_{i,j=1}^N D_j(a_{ij}(x, u) D_i u) + \frac{\varepsilon^2}{2} \sum_{i,j=1}^N D_s a_{ij}(x, u) D_i u D_j u \in (H_V(\Omega))'$$

*and $\|w_u\|_{(H_V(\Omega))'} \leqslant |df|(u)$.*

As a consequence of the previous Proposition we have the following result.

**Proposition 6.8.** *Let $u \in H_V(\Omega)$, $c \in \mathbb{R}$ and let $(u_h)$ be a sequence in $H_V(\Omega)$. Then the following facts hold:*
  (a) *if $u$ is a critical point of $f$, then $u$ is a weak solution of (62);*
  (b) *if $f$ satisfies $(CPS)_c$, then $f$ satisfies $(PS)_c$.*

Finally, we recall a suitable version of the nonsmooth Mountain–Pass Theorem.

**Proposition 6.9.** *Let us consider the class of paths (28). Assume that there exist $\varrho_\varepsilon > 0$ and $\nu_\varepsilon > 0$ such that*

(63) $$\|u\|_{H_V(\Omega)} = \varrho_\varepsilon \implies f(u) \geqslant \nu_\varepsilon.$$

*Then, if f satisfies the concrete Palais–Smale condition at level*

$$c_\varepsilon = \inf_{\gamma \in \mathscr{P}_\varepsilon} \sup_{t \in [0,1]} f(\gamma(t)),$$

*there exists a nontrivial critical point $u_\varepsilon \in H_V(\Omega)$ of f at the level $c_\varepsilon$.*

**Acknowledgments.** The author wishes to thank Marco Degiovanni for some useful discussions.

Dipartimento di Matematica e Fisica
Università Cattolica del Sacro Cuore
Via Musei 41, I–25121 Brescia, Italy.
*E-mail address*: `m.squassina@dmf.unicatt.it`